\documentclass[12pt`]{article}
\usepackage{epsfig}

\newcommand  {\Ebar} {{\mbox{\rm$\mbox{I}\!\mbox{E}$}}}
\newcommand  {\Rbar} {{\mbox{\rm$\mbox{I}\!\mbox{R}$}}}

\newcommand  {\Nbar} {{\mbox{\rm$\mbox{I}\!\mbox{N}$}}}
\newcommand {\Pbar} {{\mbox{\rm$\mbox{I}\!\mbox{P}$}}}

\newcommand{\QED}{{\hspace*{\fill}\rule{2mm}{2mm}\linebreak}}
\newsavebox{\zzzbar}
\sbox{\zzzbar}
  {\setlength{\unitlength}{0.9em}
  \begin{picture}(0.6,0.7)
  \thinlines
  \put(0,0){\line(1,0){0.6}}
  \put(0,0.75){\line(1,0){0.575}}
  \multiput(0,0)(0.0125,0.025){30}{\rule{0.3pt}{0.3pt}}
  \multiput(0.2,0)(0.0125,0.025){30}{\rule{0.3pt}{0.3pt}}
  \put(0,0.75){\line(0,-1){0.15}}
  \put(0.015,0.75){\line(0,-1){0.1}}
  \put(0.03,0.75){\line(0,-1){0.075}}
  \put(0.045,0.75){\line(0,-1){0.05}}
  \put(0.05,0.75){\line(0,-1){0.025}}
  \put(0.6,0){\line(0,1){0.15}}
  \put(0.585,0){\line(0,1){0.1}}
  \put(0.57,0){\line(0,1){0.075}}
  \put(0.555,0){\line(0,1){0.05}}
  \put(0.55,0){\line(0,1){0.025}}
  \end{picture}}
\newcommand{\Zbar}{\mathord{\!{\usebox{\zzzbar}}}}

\newtheorem{lemma}{Lemma}[section]
\newtheorem{prop}{Proposition}[section]
\newtheorem{thm}{Theorem}[section]

\newtheorem{corollary}{Corollary}[section]

\newcommand{\Z}{\Zbar}
\newcommand{\R}{\Rbar}
\newcommand{\N}{\Nbar}

\renewcommand{\L}{\ensuremath{\mathcal{L} }}
\newcommand{\E}{\Ebar}
\newcommand{\p}{\Pbar}
\newcommand{\g}{\ensuremath{(\gamma_1,\gamma_2)}}
\begin{document}

\setlength{\textheight}{21cm}

\title{{\bf On the Hydrodynamic Equilibrium of a Rod in a Lattice
    Fluid}}
\author{\bf Pablo. A. Ferrari, {\it Universidade de S\~{a}o Paulo}
  \thanks{IME USP, Caixa Postal 66281, 05315-970 - S\~{a}o Paulo, SP -
    BRAZIL - email: {\tt pablo@ime.usp.br} -
    http://www.ime.usp.br/\~{}pablo} \\
  {\bf Christian Maes {\it K.U.Leuven}}, \thanks{Onderzoeksleider FWO,
    Flanders. Instituut voor Theoretische Fysica, K.U.Leuven, Celestijnenlaan
    200D, B-3001
    Heverlee, Belgium - email:
    {\tt christian.maes@fys.kuleuven.ac.be } }\\
  {\bf Laura Ramos {\it Universidade de S\~{a}o Paulo} }\thanks{IME
    USP, Caixa Postal 66281, 05315-970
    - S\~{a}o Paulo, SP - BRAZIL - email: {\tt lramos@ime.usp.br}}, \\
  {\bf Frank Redig, {\it K.U.Leuven}}\thanks{Postdoctoraal onderzoeker FWO,
    Flanders. Instituut voor Theoretische Fysica, K.U.Leuven,
    Celestijnenlaan 200D, B-3001
    Heverlee, Belgium - email: {\tt Frank.Redig@fys.kuleuven.ac.be } }}
\maketitle

\begin{abstract}
  We model the behavior of a big (Brazil) nut in a medium of smaller
  nuts with a stochastic asymmetric simple exclusion dynamics of a
  polymer-monomer lattice system.  The polymer or `rod' can move up or
  down in an external negative field, occupying $N$ horizontal lattice
  sites where the monomers cannot enter. The monomers (at most one per
  site) or `fluid particles' are moving symmetrically in the
  horizontal plane and asymmetrically in the vertical direction, also
  with a negative field.  For a fixed position of the rod, this
  lattice fluid is in equilibrium with a vertical height profile
  reversible for the monomers' motion.  Upon `shaking' (speeding up
  the monomers) the motion of the `rod' dynamically decouples from
  that of the monomers resulting in a reversible random walk for the
  rod around an average height proportional to $\log N$.
\end{abstract}
\vspace{3mm}


\section{Introduction.}
Studying the coupled dynamics of granular matter of different
shapes and sizes is of great interest for a range of phenomena. One
example is the size segregation of particles as a result of
vibrations. A typical realization is a can with nuts; upon shaking
the larger (often taken to be Brazil) nuts rise to the top. Because
of its wide interest the phenomenon has been considered and
reconsidered and while some of the aspects are well-understood not
everything has stopped surprising.  If one asks for an analysis
starting at the microscopic level the situation is not so
satisfactory and even simple models have escaped serious
mathematical handling (cf. \cite{Leb1}, \cite{Leb2} for
further references).

In this paper we consider such a microscopic --- {\it albeit}
stochastic --- dynamics for the motion of a large particle or rod in a
lattice fluid composed of monomers. The problem of the present
paper is however not quite
similar to the canonical Brazil nuts scenario as we are interested in
the equilibrium dynamics.  In fact, as we will see, on the time scale
of the motion of the rod, the monomers are in equilibrium for a
reversible density profile.  The rod then finds its hydrodynamic
equilibrium at a vertical height where the density of the fluid is
about equal to its own density.  Going beyond equilibrium conditions,
e.g. starting from a homogeneous density for the lattice fluid, gives
rise to additional mathematical problems that we will only be touching
at the end of the paper (see Section \ref{sect:remarks}, Remark 2),
and which will be the subject of future work.\\
The result of this paper can be
classified under the heading: how to obtain a Markovian reduced
dynamics?  This problem is of course a very common one in
nonequilibrium statistical mechanics where one considers the
system composed of various types of degrees of freedom.  The
dynamics is globally defined in which the various degrees of
freedom are coupled.  In some circumstances and under some limit
procedures one then expects that some degrees of freedom of the
system effectively decouple giving rise to an autonomous (in many
cases, Markovian) dynamics for a subset of degrees of freedom.  In
our case, it is the shaking, the speeding up of the monomer
dynamics in the horizontal direction, that does the job.  In this
way, between any two moves of the polymer, the monomer
configuration has the time to relax to its reversible measure and
the polymer always sees the fluid in equilibrium.

Our main result is a mathematically rigorous proof of this
dynamical decoupling between the motion of the rod and the monomer
fluid when the monomer dynamics is (infinitely) speeded up (at
least) in the horizontal direction (orthogonal to the motion of the
polymer). In that limit of excessive horizontal shaking the reduced
dynamics of the polymer becomes that of a random walker with rates
directly given in terms of the equilibrium fluid density. When $N$
(the length of the polymer) is sufficiently big (depending on the
rates for jumping up or down) the polymer finds its most probable
height around its mean position of order $\log N$ with a variance
of order 1.  In the next section we describe the model and the
result. The third section is devoted to the proofs. The final
section contains an open problem and some additional remarks.

\section{Model and Results.}

\subsection{Model.}\label{model}
\subsubsection{Configuration.}
For convenience we put the system on the square lattice $\Z^2$.  A
point $i=(x,y)$ of the lattice has a `vertical' coordinate $y$ and a
`horizontal' coordinate $x$. We also write $i=(i_1,i_2)$ if, in the
notation of the coordinates, we want to remember the site $i$.\\ The
system contains a rigid polymer (large particle, rod) whose position
at time $t$ is denoted by $Y_t$.
For simplicity we
allow the rod to move only vertically. The horizontal coordinate is fixed
(at 0) and $Y_t$ takes values in $\Z$ (thought of as the `vertical'
axis).
The same results would hold if the polymer also jumps horizontally at
rate 1.
The polymer occupies $N\in\{ 2,3,\ldots\}$ lattice sites. If the
polymer has position $Y_t=y$, then it occupies the region
$$
A_N(y)= \{(0,y),(1,y),\ldots,(N-1,y)\}\,.
$$
  This region is forbidden for the
monomers (fluid particles).  The monomer configuration is denoted
by $\eta\in \{0,1\}^{\Z^2}$ and we use $\eta_t$ to denote the
random field of monomers at time $t$.  We have that $\eta_t(i) = 0$
if there is no monomer at site $i$ at time $t$; $\eta_t(i) = 1$ if
there is a monomer at site $i$ at time $t$.  The dynamics will
always be subject to the restriction (exclusion) that $\eta_t(i)
= 0$ for $i\in A_N(Y_t)$ (the rod acts as an obstacle for the fluid
motion). The full configuration space is denoted by $\Omega=
\{0,1\}^{\Z^2}
\times \Z$.
\\
\subsubsection{Dynamics.}
We now define the coupled dynamics for the polymer-monomers system.
All motion is via jumping to vacant sites.  There are the
horizontal jumps of the monomers (which we take symmetric and at
rate  $\gamma_1$), the vertical jumps of the monomers (asymmetric
at rate $\gamma_2$) and the vertical jumps of the rod (asymmetric
at rate 1).  The asymmetry in the vertical direction models the
presence of an external (e.g. gravitational) field acting on fluid
matter and polymer but can in general be different for monomers and
polymer. Increasing the rates $\gamma_1$ and $\gamma_2$ can be used
to simulate the greater mobility of the smaller particles upon
shaking.  We are most interested in the case where $\gamma_2 \simeq
1$ and $\gamma_1\gg 1$ (horizontal shaking).\\
\vspace{3mm}
\begin{center}
\epsfig{height=4cm, file=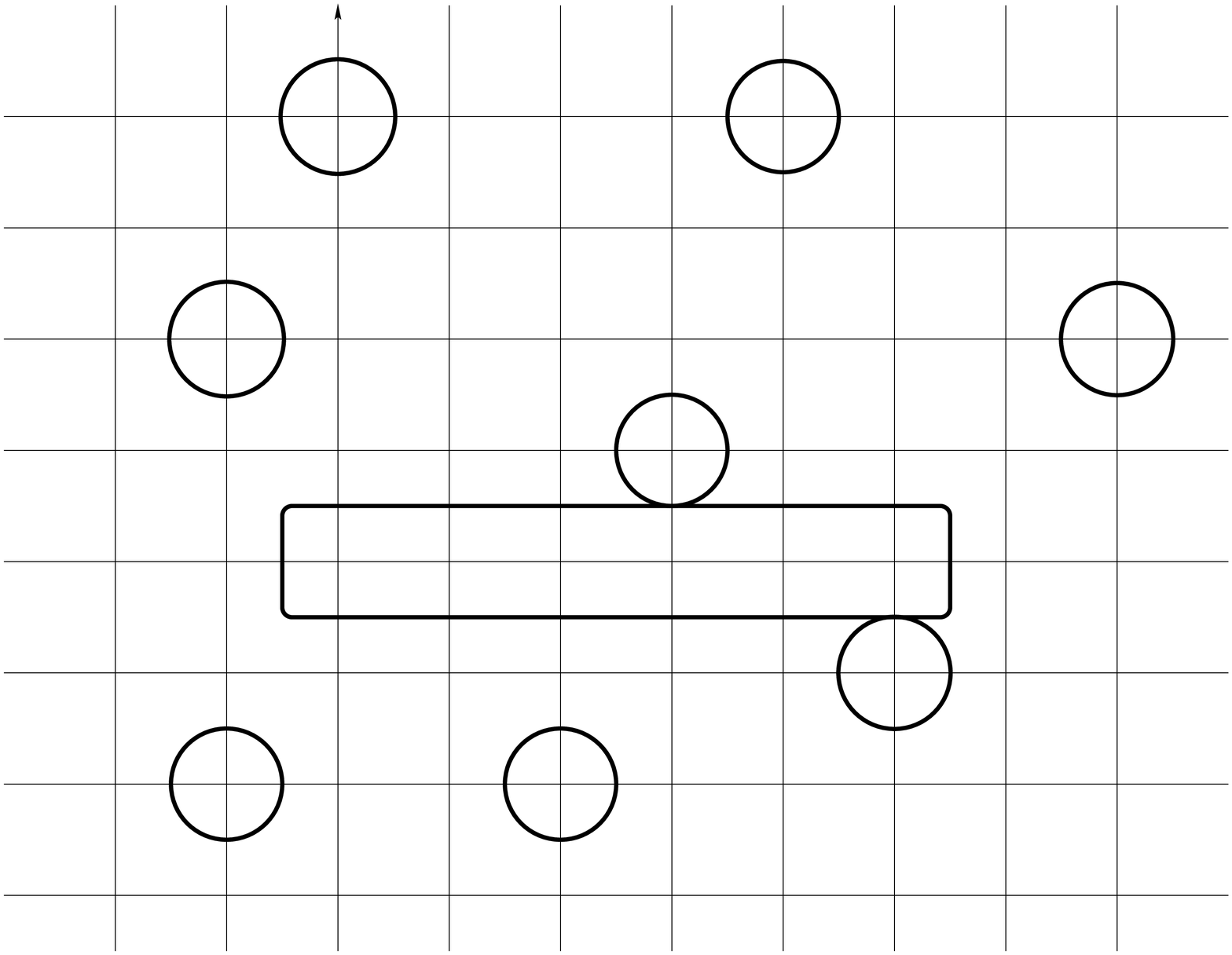}
\end{center}
\begin{center}
Fig.1: A polymer between monomers
\end{center}
\epsfig{height=4cm, file=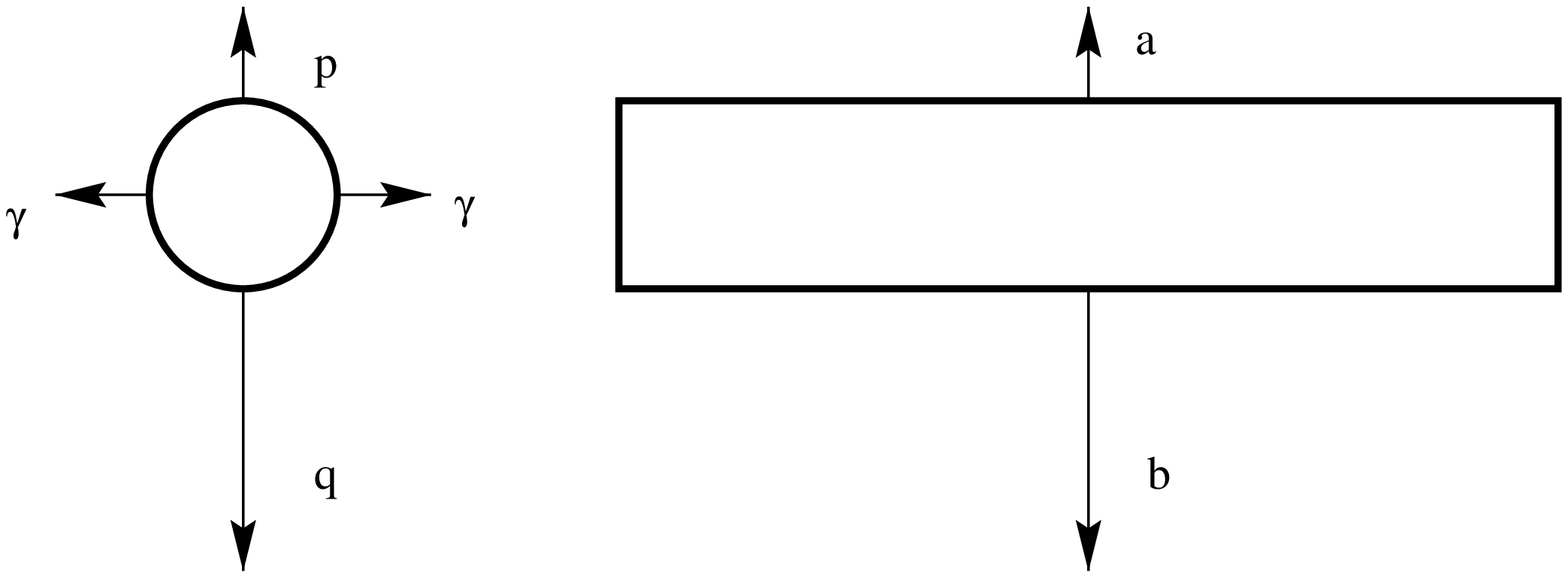}\\
\begin{center}
Fig.2: Rates of jumping for monomer and polymer where $\gamma_2=\gamma$ and
$\gamma_1=1$.
\end{center}
Here comes the formal definition of the generators of
these motions. Let $f$ be a local function on $\Omega$ (i.e., a
function that depends on the configuration in a finite region of
$\Z^2$). The first part of the generator represents horizontal monomer-jumping:
\begin{equation}\label{hori}
\L_{\rm{h}}(\eta,y) =
\frac 1{2}
\sum_{\langle ij\rangle:i_2=j_2}
I[\langle ij \rangle \cap A_N(y) = \emptyset]
[f(\eta^{i,j},y)
- f(\eta,y)]
\end{equation}
where $\eta^{i,j}(k) = \eta(i)$ if $k=j, =
\eta(j)$ if $k=i$ and $= \eta(k)$ otherwise; the summation is over
nearest neighbor pairs $\langle ij\rangle $ with the same vertical
coordinate ($i_2=j_2$). The notation $I[\cdot]$ will always stand
for the indicator function.\\ Second comes the vertical
monomer-jumping: with $p < q$,
\begin{eqnarray}\label{vert}
\L_{\rm {v}}f(\eta,y)
&=&\sum_{i} \Bigl\{p\, \eta(i)\, (1-\eta(i_1,i_2+1))\, I[(i_1,i_2 + 1) \notin
A_N(y)]\nonumber\\
&&\qquad\qquad\qquad\qquad\times[f(\eta^{i,(i_1,i_2+1)},y) - f(\eta,y)]
\nonumber\\
&&\qquad +\,\, q\, \eta(i)\, (1-\eta(i_1,i_2-1))\, I[(i_1,i_2-1)
\notin A_N(y)]\nonumber\\
&&\qquad\qquad\qquad\qquad\times[f(\eta^{i,(i_1,i_2-1)},y) - f(\eta,y)]\Bigr\}.
\end{eqnarray}
Finally, there is the polymer-jumping: with $a,b\in \R^+$,
\begin{eqnarray}
\L_{\rm {poly}} f(\eta,y) &=& a\, I[\eta(i)=0, \forall i\in A_N(y+1)]
[f(\eta,y+1)
- f(\eta,y)]\nonumber\\
&& +\; b\, I[\eta(i)=0, \forall i \in A_N(y-1)]\,
[f(\eta,y-1)
- f(\eta,y)].
\end{eqnarray}
We will then choose $p/q,a/b < 1$ to represent an external field in the
vertical direction driving all particles, big and small,
downward. E.g. in the case of a gravitational field, we could have
$p/q = \exp (-mg/kT), a/b= \exp (-Mg/kT)$, where
$m,M$ denotes the mass of a monomer, resp. polymer. \\
The formal generator $\L$ of the full dynamics consists
of three pieces:
 \begin{equation}\label{dyna}
 \L = \gamma_1 \L_{\rm {h}} + \gamma_2 \L_{\rm {v}} +
\L_{\rm {poly}},
\end{equation}
where $\gamma_1, \gamma_2 > 0$ are additional parameters governing
the rates of the monomer-jumping. Notice that $\L_{\rm{mono}} =
\gamma_1 \L_{\rm{h}} + \gamma_2 \L_{\rm{v}}$ works on the
configuration of monomers only (for fixed rod position), while
$\L_{\rm{poly}}$ works on the polymer configuration (for fixed
monomers).  The only interaction is by excluded volume. The generator
(\ref{dyna}) can be rewritten in the form:
\begin{equation}\label{rewritegen}
\L f(\eta,y )= \L_{\eta} f_\eta (y) + \L_y f_y (\eta )
\end{equation}
where
$f_\eta (\cdot) = f(\eta,\cdot )$ and $f_y (\cdot ) = f(\cdot ,y )$.

\subsubsection{Initial distribution and extra remarks.}

At time 0 (starting time) we put the polymer at the origin:
$Y_{t=0}=0$. Then fix a real parameter $\kappa$ and distribute the
monomers independently with density
\begin{equation}
\label{dens}
\rho(i) = \frac{\kappa (p/q)^{i_2}}{1 + \kappa (p/q)^{i_2}}
\end{equation}
varying in the vertical direction (constant in the horizontal
direction), conditioned on $\eta(i)=0,
\forall i \in A_N(0)$.\\
More precisely, we let $\nu_{\rho}$ denote the product measure on
$\{0,1\}^{\Z^2}$ with density
\begin{equation}
\nu_\rho(\eta(i)) = \rho(i),
\end{equation}
defined by (\ref{dens}). This measure is reversible for each of the
monomer generators process without polymer ---i.e.\/ the generators
defined by (\ref{hori}) and (\ref{vert}) but without the indicator
functions prohibiting jumps. The proof of this fact is a simple
computation. The one-dimensional analogue is well known, see
\cite{Liggett}.

For any given $y\in \Z$, we write\\
\begin{equation}
\nu_\rho^y = \nu_\rho( \cdot | \eta(i)=0, \forall i \in A_N(y)),
\end{equation}
At time 0, we put the distribution $\mu_0$ on $\Omega$ defined by
\begin{equation}
\mu_0(d\eta,y) = \delta_{y,0}\, \nu_\rho^y(d\eta),
\end{equation}
where $\delta_{y,0}$ stands for the Kronecker-delta.


{}From the initial condition described above and the dynamics defined
via (\ref{dyna}) the process $(\eta_t,Y_t)$ is generated.  The measure
at time $t \geq 0$ is denoted by $\mu_t$.  Of course this depends on
the choice of parameters $p,q,a,b,\gamma_1$ and $\gamma_2$ and we will
sometimes make this explicit in the notation.\\ A useful way to
imagine the process is by associating two exponential clocks (at rate
$a$ respectively $b$) to the polymer: one clock gives rise to the
trial times for the polymer to jump up, the other indicates the trials
for the polymer to jump down.  If, just before the trial time $\tau$, say
for jumping up, there are no monomers right above the polymer
($\eta_{\tau^-}(i)=0, \forall i \in A_N(Y_{\tau^-}+1)$), then the jump
is performed and at time $\tau$ the polymer is at height $Y_{\tau} =
Y_{\tau^-} + 1$,
otherwise it stays where it was.  Between the trial times of the
polymer, only the monomers move. The dynamics for the monomers for a
fixed position of the polymer (say at $y$) is generated by
\begin{equation}
\L_{\rm{mono}}^y f(\eta) =
\gamma_1 \L_{\rm {h}}^y + \gamma_2 \L_{\rm{v}}^y
\end{equation}
which can be read off from (\ref{hori}) and (\ref{vert}). The
associated semigroup is denoted by $S^y(t)$. Now, the important
thing where the `equilibrium' in the title of this paper refers to,
is that $\nu_\rho^y$ is a reversible measure for $S^y(t)$. This
will be proven as Lemma \ref{reversi} in Section \ref{proofs}.
\subsection{Results.}

\subsubsection{Limiting random walk.}
In the limit $\gamma_1\uparrow +\infty$ the motion of the rod
will decouple from the
monomer dynamics.  It will be a random walk.  We first
introduce this limiting rod motion.\\ For $a,b\in \R^+$ consider the
continuous time random walk on $\Z$ with generator
\begin{equation}\label{rw}
\L^{\rm RW}f(y) = a[1-\rho(y+1)]^N[f(y+1)-f(y)]+
b[1-\rho(y-1)]^N[f(y-1)-f(y)],
\end{equation}
where the density profile $\rho$ is obtained from (\ref{dens}). Remark that,
in the notation of (\ref{rewritegen}),
\begin{equation}\label{rwgen}
\L^{\rm RW}f(y) = \int \nu^y_\rho (d\eta ) \L^y_\eta f(y )
\end{equation}
$\L^{\rm RW}$
generates a continuous time
random walk $Y_t^{\rm RW}$ which we start at
$Y_{t=0}^{\rm RW}=0$ and with rate
for moving one step upward
$a[1-\rho(y+1)]^N$ and rate moving one step downward equal to
$b[1-\rho(y-1)]^N$. We fix the initial state to be $0$ for the sake of
definiteness. Our results hold for any other initial (deterministic or
random) state.
\begin{prop}
  If $a/b > (p/q)^N$, then the random walk with generator (\ref{rw})
  defined above has a unique reversible probability measure $m$ on
  $\Z$, which is given by
\begin{equation}
  \label{inv}
  m(y) = \frac{1}{Z}\,\frac{ (a/b)^y}{\left( 1+ (p/q)^y \right)^N},
\end{equation}
where $Z$ is a normalizing constant.
In particular, the random walk is positive recurrent.
\end{prop}
{\bf Proof:} Reversibility of $m(y)$ is immediate, and the condition
$a/b > (p/q)^N$ guarantees that $m(y)$ can be normalized (i.e. $Z <\infty$).
Positive recurrence follows immediately from the existence of a reversible
probability measure.\QED

Remark that the condition $a/b > (p/q)^N$ in the case of a gravitational
field just means $M/N<m$, i.e. the density of the polymer is smaller
than the density of the monomer-fluid. It is thus very natural that
in this case the polymer will drift up and will float at a height
where the fluid density is proportional to $1/N$, see (\ref{220}) and
\cite{Archi}.

In order to study some global properties of the limiting
random walk, in particular its behavior for large $N$, we
replace the discrete distribution $m(y)$ on $\Z$ by a continuous
distribution:
\begin{equation}
m(dx):= \frac{\exp (-\alpha x)}{(1+\exp{(-\beta x)})^N}
\frac{1}{Z(\alpha,\beta,N)} dx.
\end{equation}
Here
\begin{equation}\label{zet}
Z(\alpha,\beta,N) = \int_{-\infty}^\infty dx
\frac{\exp(-\alpha x)}{(1+\exp (-\beta x))^N} = \frac{1}{\beta}
\frac{ \Gamma (\alpha/\beta ) \Gamma (N-\alpha/\beta)}{\Gamma (L)},
\end{equation}
and $e^{-\alpha} = a/b$, $e^{-\beta} = p/q$.
{}From (\ref{zet}) we can calculate the cumulants of the continuous
distribution $m(dx)$: in particular
\begin{equation}\label{1moment}
\int x \,m(dx)\, =\, \frac{1}{\beta} \left( \psi (\frac{\alpha}{\beta} )-
\psi (L-\frac{\alpha}{\beta} ) \right),
\end{equation}
where $\psi (x) = \Gamma' (x)/\Gamma (x)$. Using the asymtotic
expansion
\begin{equation}
\psi (z) = \log z -\frac{1}{2z} -\frac{1}{12 z^2} + \ldots
\end{equation}
we obtain
\begin{equation}
\label{mea}
\int x\, m(dx)\, = \frac{1}{\beta}\, \log N + O(1),\ \mbox{as} \ N\rightarrow\infty
\end{equation}
and all higher order cumulants are of order 1 as $N$ tends to infinity.
The modus of $m$ (the position where $m(x)$ reaches its maximum) is
\begin{eqnarray}
\nonumber
\hbox{Mo}(m)
&=&\left(-\log (p/q)\right)^{-1}\; \log \left( \frac{\log (p/q)^N}{\log (a/b)}
- 1\right)\; \\
&\simeq& \frac{1}{\beta}\, \log N\,, \mbox{ as}\ N\rightarrow\infty.
\end{eqnarray}
\vskip 3mm

\subsubsection{Main result.}
Our main result states
\begin{thm}
\label{theo}
Let $0\le p < q<\infty$ and $a,b\in\R^+$ and consider the joint
monomer-polymer process with generator (\ref{dyna}).  For any finite
time-interval $K$, the marginal law of the polymer motion
$(Y_t^{\gamma_1}\,:\,t\in K)$ converges, as $\gamma_1\to\infty$, to
the law of the random walk $(Y_t^{\rm RW}\,:\,t\in K)$ defined by
(\ref{rw}).
\end{thm}
\subsubsection{Discussion.}
Since for $a/b > (p/q)^N$ the limiting motion is an ergodic random
walk in a countable state space, the process starting from any initial
distribution will converge to the (unique) invariant measure.  Hence,
by (\ref{mea}), the polymer will rise from the zero level to a level
at height proportional to $\frac{1}{\beta}\log N$.  If it starts in
equilibrium, then it will perform a random walk around this position.
This is exactly what we would expect from general hydrodynamics, see
\cite{Archi}.
After all, the fluid density at height $\frac{1}{\beta}\log N$ is
precisely, cf (\ref{dens}):
\begin{equation}
\label{220}
\rho (\frac{1}{\beta}\log N)\,=
\,\frac{\kappa/N}{1+\kappa/N}\,\sim\,{\kappa\over N}
\end{equation}
confirming Archimedes' law in this model of granular matter.

\section{Proofs.}\label{proofs}
\subsection{ Outline of proof}
In this section we state the main steps of the proof of
Theorem 2.1. The reader may use this section as a guideline
to the next section. 
The main idea of the proof
is that in the limit $\gamma_1\uparrow\infty$ the monomers are moving very fast
in the horizontal directions and thus can reach equilibrium in the
time between two successive jumps of the polymer.
Therefore the rate at which the polymer jumps, which is a
function of the whole monomer configuration, can be replaced by the 
expectation of that rate in the equilibrium distribution of the
monomer configuration.

As a first step (Lemma 3.1) we identify the reversible equilibrium measure
for the monomers for {\sl fixed} position of the polymer. This is (by
reversibility) the original reversible measure of the monomer gas without
polymer, conditioned on having no monomers on the lattice sites occupied
by the fixed polymer.

In a second step (Lemma 3.2-Proposition 3.1) we prove that in the
limit $\gamma_1\uparrow\infty$ any time dependent expectation of a function
$f(Y_t,\eta_t)$ of both polymer position $Y_t$ and monomer gas configuration
$\eta_t$ can be replaced by the expectation of a new function depending
only on the polymer position, and obtained from $f$ by integrating
out the $\eta$ variables over the equilibrium measure.
The main ingredients in the proof of that statement are
\begin{enumerate}
\item Discrepancies in the asymmetric exclusion process move as ``second
class particles" which are a kind of random walkers. When 
$\gamma_1\uparrow\infty$, this ``random walker" diffuses away very quickly.
\item The distribution of the monomers at any jumping time of the polymer is
absolutely continuous w.r.t. the monomer-fixed polymer equilibrium measure.
\end{enumerate}

In the first two steps we obtain convergence of the distribution of the polymer
position $Y^{\gamma_1}_t$ to the distribution of the random walk
$Y_t$. To finish our proof, we still have to prove that the {\sl whole process}
$\{ Y^{\gamma_1}_t: t\geq 0 \}$ converges
to the whole process $\{ Y_t:t\geq 0 \}$ (i.e. the distributions on trajectories
converge). This final step is made by first proving that any limiting process
is Markovian and next that there exists a limiting process (tightness).
\subsection{ Proof of Theorem 2.1}
We start this section with an easy lemma on reversible Markov processes.
\begin{lemma}\label{sile}
Let $\{ \eta_t:t\geq 0 \}$ be a Markov process on $\Omega$ with
generator $\L$ and let $\mu$ be a reversible measure for $\L$.
Suppose $A\subset \Omega$ such that $\mu (A) >0$ and such that
$1_A$ is in the domain of the generator. Consider the process with
generator
\begin{equation}\label{geno}
\L_A f = 1_A \L (1_Af)- (1_A \L 1_A)f
\end{equation}
That is, $\L_A$ corresponds to a process with ``forbidden region" $A^c$
(i.e., jumps from $A$ to $A^c$ are suppressed). Then the measure $\mu_A := \mu
(\cdot | A)$ is reversible for $\L_A$.
\end{lemma}

{\bf Proof:} Because the second term in the right hand side of (\ref{geno})
is just multiplication with the function $1_A \L(1_A)$,
it suffices to show that $\tilde{\L}_A f:=1_A \L (1_Af)$ defines
a symmetric operator on $L^2 (A,\mu_A )$.
Let $f,g$ be in the domain of $\tilde{\L}_A$. Since $d\mu_A = (1/\mu(A))1_A
\,d\mu$, we get, using the symmetry of $\L$ in $L^2(\mu )$:
\begin{eqnarray}\label{simplelemma}
\int g (\tilde{\L}_A f) d\mu_A &=& \frac{1}{\mu (A)}\int 1_A\, g
\tilde{\L}_A f\, d\mu
\nonumber\\
&=&\frac{1}{\mu (A)} \int \L(g 1_A ) 1_A f\, d\mu\nonumber\\
&=& \int \tilde{\L}_A g f \,d\mu_A.
\end{eqnarray}
\QED
Note that reversibility is crucial in the proof of this lemma. Indeed if
$\mu$ is only stationary, then we cannot conclude in general that $\mu_A$ will
be stationary for the process with generator $\L_A$. Indeed, one
easily computes
\begin{equation}
\int \L_A f d\mu_A = \frac{1}{\mu (A)} \int (1_A \L^* 1_A - 1_A \L1_A)f d\mu,
\end{equation}
i.e., $\mu_A$ will be stationary iff $1_A \L^* 1_A - 1_A \L 1_A=0$
$\mu$-a.s.
Since the profile measures are reversible for the exclusion process of the
monomers without polymer, we can apply lemma \ref{sile} for
$\mu=\nu_\rho$,
$A= \{ \eta \in \{ 0, 1\}^{\Z^2}: \sum_{z\in A(y,N)} \eta (z) \not= 0 \}$, i.e.,
those monomer configurations which are excluded when the polymer
is at vertical position $y$. This yields:
\begin{corollary}\label{reversi}
For fixed polymer position at $y \in \Z$ the measure
$\nu_\rho^y$ is reversible for the monomer dynamics with semigroup
$S^y(t)$.
\end{corollary}

\begin{lemma}\label{conv} Fix $y\in \Z$. Let $f$ be a local function
on $\{0,1\}^{\Z^2}$ which only depends on the monomer configuration
in the layers at height $y+1$ and $y-1$.  Suppose that
\[
\nu_\rho^y(f) = 0.
\]
Then, for any $t>0$,
\[
\lim_{\gamma_1 \uparrow +\infty} \| S^y(t) f\|_{L^2(\nu_\rho^y)} =
0.
\]
\end{lemma}
{\bf Proof:} Abbreviate $\mu:= \nu^y_\rho$ and consider the case
$f_y (\eta )= 1_A (\eta ) - \mu (A)$ for a set $A$ in the space of
configurations depending only on a finite number of coordinates in
labels $y-1$ and $y+1$. The extension to general local $f$
is straightforward. Denote
\begin{equation}\label{D_A}
D_A := \{ x\in \Z^2 : 1_A (\eta ) \not= 1_A (\eta^x )
\ \mbox{for some} \ \eta \},
\end{equation}
the dependence set of $A$. By reversibility:
\begin{eqnarray}\label{tralala}
\int  (S^y (t ) f_y )^2 d\mu &=& \int (S^y (2t) f_y) f_y d\mu.
\nonumber \\
&=& \mu (A) \left( \E^y_{\mu(\cdot |A)} (1_A (\eta_t )) -
\E^y_\mu (1_A (\eta_t ))\right).
\end{eqnarray}
To compute the difference of the expectations in the above expression
we realize the processes with initial configurations $\eta$ and
$\zeta$ in the same probability space (coupling).

To construct this coupling we first associate two Poisson clocks to
each site of $\Z$ with parameters $\gamma_2p$ and $\gamma_2q$
respectively and use them to decide the times of the vertical
attempted jumps. A jump from $(i_1,i_2)$ to $(i_1,i_2+1)$ is performed
at time $t$ if an event of the Poisson process of rate $p$ occurs at
that time, a particle is present at $(i_1,i_2)$ and no particle is
present at $(i_1,i_2+1)$ at time $t-$. Similarly, a jump from
$(i_1,i_2)$ to $(i_1,i_2-1)$ is performed at time $t$ if an event of
the Poisson process of rate $q$ occurs at that time, a particle is
present at $(i_1,i_2)$ and no particle is present at $(i_1,i_2-1)$ at
time $t-$. Jumps either to or from sites occupied by the rod are
suppressed. This takes care of the vertical jumps. See Ferrari (1992)
for details of this construction. For the horizontal jumps we
associate Poisson clocks with rate $\gamma_1$ to pairs of horizontal
nearest-neighbor sites. When the clock associated with sites
$(i_1,i_2)$ and $(i_1+1,i_2)$ rings, the contents of those sites are
interchanged.  Also here, if at least one of the sites is occupied by
the rod, the jump is suppressed. The horizontal motion is also called
\emph{stirring} process. See Arratia (1986) for details of this
construction.  More rigourosly, let $(N_t(i,j)\,:\, i=(i_1,i_2)\in \Z,
j=(i_1+1,i_2))$, $(N^+_t(i)\,:\, i\in \Z)$ and $(N^-_t(i)\,:\, i\in
\Z)$ three independent families of independent Poisson processes of
rates $\gamma_1/2$, $p\gamma_2$ and $q\gamma_2$, respectively ---the
Poisson clocks. Use the notation ${\rm d}N_t(\cdot)=1$ if there is an event
of the Poisson process $(N_t(\cdot))$ at time $t$, otherwise it is
zero. The motion is defined by
\begin{eqnarray}
  \label{g1}
 {\rm d}f(\eta_t) &=& \sum_{\langle ij\rangle:i_2=j_2} {\rm d}N_t(i,j)\,
I[\langle ij \rangle \cap A_N(y) = \emptyset]\,
[f(\eta^{i,j},y)
- f(\eta,y)]\nonumber\\
&&+\,\, \sum_{i} \Bigl\{{\rm d}N^+_t(i)\, \eta(i)\, (1-\eta(i_1,i_2+1))\,
I[(i_1,i_2 + 1) \notin A_N(y)]\nonumber\\
&&\qquad\times\,[f(\eta^{i,(i_1,i_2+1)},y) - f(\eta,y)]
\nonumber\\
&& +\,\, {\rm d}N^-_t(i)\, \eta(i)\, (1-\eta(i_1,i_2-1))\,
I[(i_1,i_2-1) \notin
A_N(y)]\nonumber\\
&&\qquad\times\,[f(\eta^{i,(i_1,i_2-1)},y) - f(\eta,y)]\Bigr\}.
\end{eqnarray}
Standard arguments, see for instance Durrett (1993) show that
(\ref{g1}) defines a process $\eta_t = \Phi (\eta_0;N[0,t])$, with
initial configuration $\eta_0$, where $\Phi$ is the function induced
by (\ref{g1}) and $N[0,t]\,:=\,(N_s(\cdot,\cdot),\,
N^+_s(\cdot),\,N^-_s(\cdot)\,:\, 0\le s\le t)$; furthermore it is
immediate to see that $\eta_t$ has generator $\L_y$.  Given two
initial configurations $\eta$ and $\zeta$, the coupling of their
evolutions is constructed using the same Poisson processes: define
$$
(\eta_t,\zeta_t):= (\Phi (\eta;N[0,t]),\Phi (\zeta;N[0,t])).
$$
Let $\E^y_{(\eta,\zeta)}$ denote expectation in the coupling
starting with $(\eta,\zeta)$. We need also to couple the initial
configurations. Let $\widetilde\mu_A$ be the law of a pair of
configurations $(\eta,\zeta)$ with marginal distributions $\mu_A$ and
$\mu$ and such that $\eta(x)=\zeta(x)$ for all $x\in \Z^2\setminus
D_{A}$. It is possible to construct a measure with these properties
because $\mu$ is a product measure. We then have
\begin{equation}
  \label{co2}
  \E^y_{\mu(\cdot |A)} (1_A (\eta_t )) -
\E^y_\mu (1_A (\eta_t )) = \int \widetilde\mu_A(d(\eta,\zeta))\,
\E_{(\eta,\zeta)} [1_A (\eta_t ) - 1_A (\zeta_t )]
\end{equation}

The number of initial discrepancies is finite, that is, $\sum_{x} I
(\eta (x) \not= \zeta (x) )\,\le\, |D_A|\,<\, \infty$. At each site
$x$ of $\Z$ we have one of three possibilities:
$(\eta(x)-\zeta(x))=0$, no discrepancies; $(\eta(x)-\zeta(x))^+>0$,
positive discrepancies; or $(\eta(x)-\zeta(x))^->0$, negative
discrepancies. Following the evolution of the particles and the
discrepancies we notice that if a positive
discrepancy jumps over a negative one, then both discrepancies
collide, giving place to a coupled particle and a hole; if a coupled
particle attempts to jump to a discrepancy, the jumps occur and then
the discrepancy must jump to the site previously occupied by the
coupled particle. These two behaviors only occur when vertical jumps
are involved. In the horizontal jumps, discrepancies and coupled
particles just interchange positions according to the Poisson
horizontal (stirring) clocks.

We say that there is a \emph{first class particles}
at site $i$ at time $t$ when
$\xi_t(i)=\eta_t(i)\zeta_t(i)=1$, a \emph{positive second class particles}
when $(\eta\zeta)_t(i) = \eta_t(i)-\zeta_t(i)=1$ and
a \emph{negative second class particles} when
$(\zeta\eta)_t(i)=\zeta_t(i)-\eta_t(i)=1$. The first class particles
occupy initially those sites $i$ occupied by both $\eta$ and
$\zeta$. Locally in time, the motion of the
first class  particles is the one given by generator $\L^y$ but
superposed to it there is a pure birth process of first class
particles: with
rate
$$
p\,(\eta\zeta)_t(i_1,i_2-1)\,(\zeta\eta)_t(i_1,i_2)
$$
the second class particles at $(i_1,i_2-1)$ and $(i_1,i_2)$ annihilate
each other and a first class particle appears at
$(i_1,i_2)$ and an empty site appears at $(i_1,i_2-1)$. Similarly, at rate
$$
q\,(\eta\zeta)_t(i_1,i_2+1)\,(\zeta\eta)_t(i_1,i_2)
$$
the second class particles at $(i_1,i_2+1)$ and $(i_1,i_2)$ annihilate
each other and a first class particle appears at
$(i_1,i_2)$ and an empty site appears at $(i_1,i_2-1)$.

The marginal distribution of a second class particle between two
vertical jumps (or between a jump and an annihilation) corresponds to
the law of a nearest neighbor symmetric random walk ---with reflection
at the rod when at level $y$--- in the horizontal direction. In the
vertical direction the motion is not Markovian ---it depends on the
configuration of the first and second class particles at the instants
of attempted jumps--- and either there is an annihilation as described
above or the second class particles just change horizontal line. For
instance, at time $t$, jumps of a $(\eta\zeta)$ second class particle
{}from site $(i_1,i_2)$ to site $(i_1,i_2+1)$ occur with rate
$$
p\,(\eta\zeta)_t(i_1,i_2)\,(1-\xi_t(i_1,i_2+1)) + q\,
(\eta\zeta)_t(i_1,i_2)\,\xi_t(i_1,i_2+1)
$$
and similarly for the other cases. The first term corresponds to the
jump over an empty site and the second one to the interchange of
positions with a first class particle.

This coupling has the property
\begin{equation}
\p^y_{(\eta,\zeta)} \Bigl\{ \sum_{x} I (\eta (x) \not= \zeta (x) )
\ge \sum_x I (\eta_t (x) \not= \zeta_t (x) ) \Bigr\}=1,
\end{equation}
i.e., the number of discrepancies cannot increase.

Since by construction the discrepancies between $\eta$ and $\zeta$ are
all located at $D_A$, we have the estimate
\begin{equation}
\label{esti1}
\E^y_{(\eta,\zeta)} [ 1_A(\eta_t ) -1_A (\zeta_t) ]
\leq \sum_{i\in D_A}\sum_{z\in D_A} \p ( X^{i}(t) = z),
\end{equation}
where $X^{i}(t)$ is the position of a second class particle
initially at $i$. If at site $i$ there were no discrepancy we use
the convention $X^i_t\not\in \Z^2$ (and hence $\neq z$, for all
$z\in D_A$). If particles $i$ and $j$ were discrepancies of
different sign and collided before time $t$, we also set
$X^i_t,X^j_t\not\in \Z^2$. The process $X^i_t$ has rate $\gamma_1$
to move symmetrically in the horizontal direction. If the rod were
not present, we could dominate $\p ( X^{i}(t) = y)$ by
$\p((X^i(t))_1 = y_1)$, where $(X^i(t))_1$ is the first coordinate
of the walk. Since without the rod the first coordinate makes just
a symmetric random walk at rate $\gamma_1$, that probability would
be dominated by $\gamma_1^{-1/2}$ times a constant. But with the
rod we have to work a bit more. The process $X^i(t)$ has rate at
most $\gamma_2(p+q)$ to move in the vertical direction. This
implies that the time elapsed between the last vertical jump and
$t$ is dominated by the minimum between an exponential time of rate
$\gamma_2(p+q)$ and $t$. With this in hand it is not difficult to
prove that also in this case $\p ( X^{i}(t) = y)$ is bounded above
by $\gamma_1^{-1/2}$ times a(nother) constant. Here we use that
$\gamma_2$ remains bounded when $\gamma_1$ goes to infinity. We
conclude that for any pair $i$, $y$ in $D_A$:
\begin{equation}
  \label{wcn}
 \lim_{\gamma_1\uparrow\infty} \p ( X^{i}(t) = y)\,=\,0
\end{equation}

Therefore we conclude, combining
(\ref{tralala}), (\ref{co2} ), (\ref{esti1} ) and (\ref{wcn}) and the
fact that $D_A$ is a finite set:
\begin{eqnarray}
\lim_{\gamma_1\uparrow\infty} \int d\mu\, f_y\, S^y (2t ) f_y \;=\;0 .
\end{eqnarray}
\QED

{\bf Remark:} We postpone until Section \ref{sect:remarks}, Remark 3,
an alternative more general proof of Lemma \ref{conv} which works
equally well for a broader class of exclusion dynamics (e.g. with
speed change) provided the projection of
the invariant measure on horizontal layers
is ergodic for the horizontal dynamics.

We now prove an intermediate result which is
important for the proof of Theorem~2.1.
\begin{prop}\label{interprop}
  Let $f_y$ be a function depending only on the configuration values
  at the $N$ sites of $A_N(y-1)$ or $A_N(y+1)$. Then we have for all
  $t>0$:
\begin{equation}
\limsup_{\gamma_1\uparrow\infty}
\E^{(\gamma_1,\gamma_2)}_{\nu^0_\rho \times\delta_0}
\left( f_{Y_t} (\eta_t ) - \int f_{Y_t} (\eta )\, \nu^{Y_t}_{\rho} (d\eta )
\right) =0.
\end{equation}
\end{prop}
{\bf Proof:} We first want to condition on a sequence $\underline
T^\epsilon:=(T^{\epsilon_1}_1,\ldots, T^{\epsilon_n}_n)$ of marked
trial jumps before $t$. Here
$\epsilon\in \{ -1 , +1 \}$ is the mark of the jump: $+1$ for up, $-1$
for down. Next we consider $\alpha_1,\ldots, \alpha_n \in \{0,1 \}$
with interpretation $\alpha_i=1$ if $i$-th marked trial jump succeeds,
$\alpha_i =0$ if not. Given
$(T^{\epsilon_1}_1,\ldots,T^{\epsilon_n}_n)$ and
$\alpha:=(\alpha_1,\ldots,\alpha_n)$, we define
\begin{equation}
Y^\alpha_k = \sum_{j=1}^k \epsilon_j \alpha_j.
\end{equation}
This corresponds to the position of the polymer at time
$T^{\epsilon_k}_k$, given succeeded and failed jumps
$(\alpha_1,\ldots,\alpha_k)$. Finally we denote by
$V^{\alpha,\epsilon}_p:=V_{\underline T^\epsilon,\alpha_p}^\alpha (\eta )$
the event that the polymer in $Y^\alpha_{p-1}$ can (for $\alpha_p=1$)
or cannot (for $\alpha_p=0$) perform the jump to $Y^\alpha_{p-1}
+\epsilon_p$.  With this notation, we can write
\begin{eqnarray}
\lefteqn{\hskip-5mm\E^{(\gamma_1,\gamma_2)}_{\nu^0_\rho\times\delta_0}
\left(   f_{Y_t} (\eta_t ) - \int f(\eta ) \nu^{Y_t}_\rho (d\eta )
\,\Bigm|\, T^{\epsilon_1}_1,\ldots, T^{\epsilon_n}_n \,;\,
T^{\epsilon_n}_n < t < T^{\epsilon_{n+1}}_{n+1}
\right)}\nonumber\\
&=& \sum_{\alpha\in\{ 0,1 \}^{ \{1,\ldots,n\}}}
\p^{(\gamma_1,\gamma_2)}_{\nu^0_\rho\times\delta_0}
(\alpha) \label{step1}\\
&&\qquad\times\int d{\mu^{(\gamma_1,\gamma_2)}_{\alpha,Y^\alpha_n}}
S_{(\gamma_1,\gamma_2)}^{Y^\alpha_n}( t - T^{\epsilon_n}_n )
\left( f_{Y^\alpha_n} -
\int f_{Y^\alpha_n} (\eta ) \nu^{Y^\alpha_n}_\rho (d\eta ) \right).
\nonumber
\end{eqnarray}
Here $\mu^{\g}_{\alpha,Y^\alpha_n}$ denotes the monomer distribution
at time $s=(T^{\epsilon_n}_n)^+$, given the successes
$(\alpha_1,\ldots,\alpha_n)$, and
$\p^{\g}_{\nu^0_\rho\times\delta_0}(\alpha )$ denotes the probability
of the sequence of succeeded and failed jumps prescribed by $\alpha$
at the times $\underline T^\epsilon$.
The crucial thing to realize at this point is that the
probability measure $\mu^{\g}_{\alpha,Y^\alpha_n}$ is absolutely
continuous with respect to the conditioned Bernoulli measure
$\nu^{Y^\alpha_n}_\rho$. In Lemma (\ref{esti}) below we shall give a
uniform bound on the density
\begin{equation}
\Psi^{\g}_{\alpha,n} := \frac{d\mu^{\g}_{\alpha,Y^\alpha_n}}{d\nu^{Y^\alpha_n}_\rho}.
\end{equation}
By dominated convergence, the proof of the proposition is reduced to showing
that for any $\alpha\in\{ 0,1 \}^ { \{1,\ldots, n\}}$ and any $\delta > 0$:
\begin{equation}\label{afsch1}
\lim_{\gamma_1\uparrow\infty}
\int d\mu^{\g}_{\alpha,Y^\alpha_n} S^{Y^\alpha_n}_{\g} (\delta )
\left( f_{Y^\alpha_n} - \int f(\eta ) \nu^{Y^\alpha_n}_\rho (d\eta )
\right)= 0.
\end{equation}
The expression inside the limit in the left hand side of
(\ref{afsch1}) is bounded by
\begin{equation}
\| \Psi^{\g}_{\alpha,n} \|_\infty \;\Bigl\| S^{Y^\alpha_n}_{\g} (\delta )
\,\Bigl( f - \int f(\eta )\,\nu^{Y^\alpha_n}_\rho (d\eta )\Bigr)
\Bigr\|_{L^2 \left(\nu^{Y^\alpha_n}_\rho \right)}.
\end{equation}
Therefore, (\ref{afsch1}) is a consequence of Lemma \ref{conv}
and the following estimate on the density $\Psi^{\g}_{\alpha,n}$.
\begin{lemma}\label{esti}
Put $c(\rho,x ):= [\rho (x+1) \wedge \rho (x-1)\wedge (1-\rho (x))]^{-N}$.
For any $\alpha\in \{ 0, 1 \}^{\N}$
and for any $n\in\N$,
we have the estimate:
\begin{equation}\label{esti3}
\limsup_{\gamma_1\uparrow\infty} \| \Psi^{\g}_{\alpha,n} \|
\;\leq\; \prod_{p=0}^{n-1} c(\rho, Y^\alpha_p )
\end{equation}
\end{lemma}
{\bf Proof:} We fix $\alpha$ and proceed by induction in $n$. First put $n=1$.
By stationarity of $\nu^0_\rho$ under the evolution $S^0_{\g}$, we have
\begin{equation}
\mu^{\g}_{\alpha_1,Y^{\alpha_1}_1} = \nu^0_\rho [\, \cdot\, | V^{\alpha}_1 ].
\end{equation}
First consider $\alpha_1 =1$, i.e., the jump succeeds. Denote
$V(x)$ the event that the set
$A_N (x)$ contains no monomers. Then we can write:
\begin{equation}
\int f(\eta )\nu^0_\rho  [d\eta|V^{\alpha_1}_1]
= \int d\nu^{Y^{\alpha_1}_1}_\rho [ f I(V(0))] \frac{\nu_\rho (V(Y^{\alpha_1}_1))}
{\nu_\rho (V (0))}.
\end{equation}
Hence, we conclude
\begin{equation}
\Psi^{\g}_{\alpha_1} =  I ( V(0) ) \frac{\nu_\rho (V(Y^{\alpha_1}_1 ))}
{\nu_\rho ( V(0) )}.
\end{equation}
And we can estimate
\begin{equation}
\| \Psi^{\g}_{\alpha_1} \| \leq \frac{1}{\nu_\rho ( V(0) )} \leq c(\rho, 0 ).
\end{equation}
Next consider $\alpha_1 = 0$, i.e., the jump fails (and thus
$Y^{\alpha_1}_1 =0$). We write
\begin{equation}
\int f(\eta )\nu^0_\rho [d\eta| V^{\alpha_1}_1] = \int\nu^{Y^{\alpha_1}_1}_\rho (d\eta)
[f I (V^{\alpha_1}_1 (\eta ))] \frac{1}{\nu_\rho^{Y^{\alpha_1}_1}
(V^{\alpha_1}_1)}.
\end{equation}
Hence,
\begin{equation}
\Psi^{\g}_{\alpha_1}= \frac{ I [ V^{\alpha_1}_1 ]}{\nu^0_\rho ( V^{\alpha_1}_1
)}
\end{equation}
So also in that case we have the estimate
\begin{equation}
\| \Psi^{\g}_{\alpha_1} \| \leq \frac{1}{\nu^0_\rho (V^{\alpha_1}_1)}\leq
c(\rho,0 ).
\end{equation}
This proves the claim for $n=1$. Suppose the claim is true for
$n=1,\ldots, p-1$. Put $\alpha_p = 1$, the case $\alpha_p =0$ can be
treated analogously. In order to simplify the notation, we make some further
abbreviations:
\begin{enumerate}
\item $\mu^{\g}_{\alpha,Y^{\alpha}_p}:= \mu_p$.
\item $\nu^{Y^{\alpha_1,\ldots,\alpha_p}_p}_\rho:= \nu^p_\rho$
\item $\Psi^{\g}_{\alpha_1,\ldots,\alpha_p}:= \Psi^{\g}_p$
\item $S^{Y^{\alpha_1,\ldots,\alpha_p}_p}_{\g} (t):=S^{\g}_p (t)$
\item $T^{\epsilon_p}_p - T^{\epsilon_{p-1}}_{p-1} := \tau_p$
\item $V^{\alpha_1,\ldots,\alpha_p}_p:=V_p$
\end{enumerate}
We compute $\Psi^{\g}_p$:
\begin{eqnarray}\label{esti2}
\mu_p (f) &=& \left( \mu_{p-1} S^{\g}_{p-1} (\tau_p ) \right)[f|V_p]\nonumber\\
&=&
\frac{\int d\mu_{p-1}\  S^{\g}_{p-1} (\tau_p ) (f 1_{V_p})}
{\int d\mu_{p-1}\  S^{\g}_{p-1} (\tau_p ) (1_{V_p} )}
\nonumber\\
&=&
\frac {\int d\nu^{p-1}_\rho \left( S^{\g}_{p-1} (\tau_p ) (\Psi^{\g}_{p-1} )\
f 1_{V_p} \right)}
{\int d\nu^{p-1}_\rho \left( S^{\g}_{p-1} (\tau_p ) (\Psi^{\g}_{p-1})\ 1_{V_p}
\right)}
\nonumber\\
&=&
\frac {\int d\nu^p_\rho \left( 1_{V_{p-1}} f\
S^{\g}_{p-1} (\tau_p ) (\Psi^{\g}_{p-1})\right)}
{\int d\nu^p_\rho \left( 1_{V_{p-1}} S^{\g}_{p-1} (\tau_p ) (\Psi^{\g}_{p-1})
\right)},
\end{eqnarray}
where in the third step we used reversibility of $\nu^{p-1}_\rho$.
{}From (\ref{esti2}) we read off the density:
\begin{equation}\label{step2}
\Psi^{\g}_p =  \frac {1_{V_{p-1}}\ S^{\g}_{p-1} (\tau_p ) (\Psi^{\g}_{p-1} )}
{\int d\nu^p_\rho \left(1_{V_{p-1}}\ S^{\g}_{p-1} (\tau_p ) \Psi^{\g}_{p-1}
\right)}.
\end{equation}
We first estimate the nominator of the rhs of (\ref{step2}):
\begin{eqnarray}\label{step3}
&&\int d\nu^p_\rho \ 1_{V_{p-1}} \ S^{\g}_{p-1} (\tau_p ) (\Psi^{\g}_{p-1} )
\nonumber\\
&=&
\int d\nu^{p-1}_\rho \left( 1_{V_p} S^{\g}_{p-1} (\tau_p ) (\Psi^{\g}_{p-1} )
\right) \frac {\nu_\rho ( V_{p-1} )}{\nu_\rho (V_p )}\nonumber\\
&=& \int d\nu^{p-1}_\rho \left( S^{\g}_{p-1} (\tau_p ) (1_{V_p} )
\  \Psi^{\g}_{p-1}\right)
\frac {\nu_\rho ( V_{p-1} )}{\nu_\rho (V_p )}\nonumber\\
&\geq & \nu_\rho (V_{p-1} )\nonumber\\
&-&\frac {\nu_\rho ( V_{p-1} )}{\nu_\rho (V_p )}
\| \Psi^{\g}_{p-1} \|_{\infty} \| S_{p-1}^{\g} (\tau_p )
[1_{V_p} -\nu^{p-1}_\rho (V_p )] \|_{L^2 (\nu^{p-1}_\rho )}\nonumber\\
&\geq & \frac{1}{c(\rho,Y^{\alpha}_{p-1} )} - o(\gamma_1 ),
\end{eqnarray}
where $o(\gamma_1 )$ tends to zero as $\gamma_1 \uparrow \infty$ by
Lemma \ref{conv}.
By the induction hypothesis, we obtain from (\ref{step2}), (\ref{step3}):
\begin{eqnarray}
\limsup_{\gamma_1\uparrow\infty } \| \Psi^{\g}_p \|_{\infty}
&\leq & \limsup_{\gamma_1 \uparrow\infty} \|\Psi^{\g}_{p-1} \|_{\infty}
c(\rho, Y^{\alpha}_{p-1} )
\nonumber\\
&\leq & \prod_{k=0}^{p-1} c(\rho,Y^{\alpha}_k ).
\end{eqnarray}
This finishes the proof of Lemma
\ref{esti} and Proposition \ref{interprop}.\QED
As a first application we obtain convergence of the
one-point marginales of the processes $\{ Y^{\g}_t:t\geq 0\}$. For
$f:\Z\rightarrow\R$ a bounded function, we have, using the notation
of (\ref{rewritegen}).
\begin{equation}
\E^{\g}_{\nu^0_\rho\times\delta_0}\left( f(Y_t)-f(Y_0)
-\int_0^t ds\  (\L_{\eta_s} f)(Y_s ) \right)=0.
\end{equation}
By Proposition \ref{interprop} we obtain in the limit $\gamma_1\uparrow\infty$:
\begin{eqnarray}\label{onedimmarg}
&&\lim_{\gamma_1\uparrow\infty}
\E^{\g}_{\nu^0_\rho\times\delta_0} \left( f(Y_t) -f(Y_0)
-\int_0^t ds\  [(\L_\eta )\nu^{Y_s}_\rho (d\eta )]f(Y_s)\right)
\nonumber\\
&=& \lim_{\gamma_1\uparrow\infty} \E^{\g}_{\nu^0_\rho\times\delta_0}
\left( f(Y_t)-f(Y_0)-\int_0^t ds\ (\L^{\rm RW}f)(Y_s) \right)=0.
\end{eqnarray}
This implies in particular that
\begin{equation}
\lim_{\gamma_1\uparrow\infty} \E^{\g}_{\nu^x_\rho\times\delta_x} f(Y_t)
= \E^{\rm RW}_x f(Y^{\rm RW}_t).
\end{equation}
In order to prove that the {\it processes} $\{ Y^{\g}_t:t\geq 0 \}$ converge
weakly in the Skorohod space of trajectories to the random walk
$\{ Y^{\rm RW}_t:t\geq 0\}$, i.e. the content of Theorem 2.1, it suffices
to show that the process $\{ Y^{\g}_t: t\geq 0 \}$ is asymptotically Markovian.
Indeed, then it is uniquely determined by its single time distributions 
which are
those of the random walk $\{ Y^{\rm RW}_t: t\geq 0 \}$.
More precisely it is sufficient to prove the following lemma:
\begin{lemma} Let $\{ \mathcal{F}_t:t\geq 0 \}$ denote the $\sigma$-field
generated by $\{ (\eta_s,Y_s ):0\leq s\leq t \}$. We have
\begin{equation}\label{bijna}
\lim_{\gamma_1\uparrow\infty} \E^{\g}_{\nu^0_\rho\times\delta_0}
\left|
\E^{\g}_{\nu^0_\rho\times\delta_0} \left( f(Y_t)|\mathcal{F}_s \right)
-
\sum_{y} p^{\rm RW}_{t-s} (Y_s,y) f(y) \right| =0
\end{equation}
\end{lemma}
{\bf Proof:} By the Markov property of the process $\{ (\eta_t, Y_t ): t\geq
0\}$,
\begin{eqnarray}\label{markovprop}
\E^{\g}_{\nu^0_\rho\times\delta_0} \left( f(Y_t) | \mathcal{F}_s \right)
&=& \E^{\g}_{\eta_s\times\delta_{Y_s}} (f ( Y_{t-s} ))\nonumber\\
&=& f(Y_s ) + \E^{\g}_{\eta_s\times\delta_{Y_s}}
\int_0^{t-s} \L_{\eta_r} f ( Y_r ) dr.
\end{eqnarray}
Therefore, it suffices to show that
\begin{equation}\label{sufficient}
\lim_{\gamma_1\uparrow\infty} \E^{\g}_{\eta_s\times\delta_{Y_s}}
\left(\int_0^{t-s} \L_{\eta_r} f ( Y_r ) dr
-\int_0^{t-s} dr \int \L_{\eta} f(Y_r ) \nu^{Y_r}_\rho (d\eta ) \right)=0
\end{equation}
Since the trial jumps of the polymer are on the event times of a Poisson
process with rate independent of $\g$, we can write
\begin{equation}
\int_0^{t-s} \L_{\eta_r} f (Y_r ) dr
= \int_0^{t-s} \frac{1}{\epsilon} \int_r^{r+\epsilon}
\L_{\eta_{r'}} f(Y_{r'}) dr' + o(\epsilon ),
\end{equation}
where $o(\epsilon )$ goes to zero in $L^2
(\p^{\g}_{\nu^0_\rho\times\delta_0})$, uniformly in $\g$, when
$\epsilon$ tends to zero.  Therefore,
it is sufficient to show that
\begin{equation}
\lim_{\gamma_1\uparrow\infty} \left(\E^{\g}_{\nu^0_\rho\times\delta_0}
\E^{\g}_{\eta_s\times\delta_{Y_s}}\left| \frac{1}{\epsilon}
\int_0^\epsilon f_{Y_r} (\eta_r ) dr -
\int \nu^{Y_s}_\rho (d\eta ) f_{Y_s} (\eta)\right|\right)=0
\end{equation}
Following the same strategy as in the proof of Proposition \ref{interprop},
i.e., by estimates on the density of the monomer distribution with respect to
the appropriate conditioned Bernoulli measure, this reduces to showing that
for any $\epsilon >0$, for any $y\in \Z$ and for $f_y$ depending on
layer $y+1$ or $y-1$:
\begin{equation}\label{ok}
\lim_{\gamma_1\uparrow\infty} \E^{\g,y}_{\nu^y_\rho}
\left( \frac{1}{\epsilon} \int_0^\epsilon ds\  f_y (\eta_s )
- \int \nu^y_\rho (d\eta ) f_y (\eta ) \right)^2 =0.
\end{equation}
Putting $\tilde{f}_y := f_y - \nu^y_\rho (f)$, the expression inside
the limit in the left hand side of (\ref{ok}) can be rewritten as
\begin{eqnarray}
&&\int \nu^y_\rho (d\eta ) \frac{1}{\epsilon^2} \int_0^\epsilon ds
\int_0^\epsilon dr \tilde{f}_y S_{\g}^y (|r-s| ) \tilde{f}_y
\nonumber\\
&\leq & \frac{1}{\epsilon^2} \int_0^\epsilon ds\int_0^\epsilon dr
\|\tilde{f}_y \|_{L^2} \| S_{\g}^y (|r-s| ) \tilde{f}_y \|_{L^2}.
\end{eqnarray}
Hence we obtain (\ref{ok}) as an application of Lemma \ref{conv}.\QED
Arrived at this point, we know that any weak limit point of the
processes $\{ Y^{\g}_t :t\geq 0 \}$ equals in distribution
the random walk $\{ Y^{\rm RW}_t :t\geq 0 \}$. Hence, to finish the proof
of Theorem \ref{theo}, it is sufficient to see that such a weak limit point
actually exists. This is an easy task:
\begin{lemma}\label{tight}
The sequence of processes $\{ Y^{\g}_t:t\in [0,T], \}_{\g}$ is tight.
\end{lemma}
{\bf Proof:}
Since the number of jumps the polymer makes in [0,T] is bounded by
a mean one Poisson process, we have
\begin{equation}
\p (\sup_{0\leq t\leq T} |Y^{\g}_s| \geq M ) \leq \frac{2T}{M},
\end{equation}
and also
\begin{equation}
\lim_{\delta\downarrow 0}
\p (\sup_{s,t\in [0,T], |s-t|\leq \delta} |Y^{\g}_s-Y^{\g}_t|>\epsilon)
=0.
\end{equation}
This proves tightness (cf. Theorem 1.3 p.51 of \cite{Landim}).\QED
\section{Additional remarks}\label{sect:remarks}
\noindent {\bf Remark 1:}
What happens when the system is out of equilibrium? For instance,
start the monomers in a homogeneous product measure.  When the
density is constant and equal to $\rho\in[0,1]$ (no $p,q,i_2-$
dependence in (\ref{dens}), the measure $\nu_\rho^y$ is no longer
invariant for the monomer dynamics $S^y(t)$ at fixed rod position
$y$. However in the limit $\gamma_1\to\infty$ the polymer will
perform a continuous time random walk with rates $a(1-\rho)^N$ and
$b(1-\rho)^N$ for up and down jumps respectively.  Significant
corrections in the case $\gamma_1 < \infty$ can be expected, cf.
\cite{Leb2}.\\
Another problem is obtained if we start the
monomers from a sharp density profile.  That is, above the polymer the
fluid density is constant $\rho_1$ and under the polymer the density
is also homogeneous equal to $\rho_2$. In this case the vertical
density will follow a discrete space noiseless Burgers equation:
$\rho(i,t)\in[0,1]$, $t\in \R$, $i\in \Z$
\begin{eqnarray}
  {\partial \rho(i,t) \over \partial t}
&=& - p \,  \rho(i,t)\,(1-\rho(i+1,t))\,
-\,q\, \rho(i,t)\,(1-\rho(i-1,t))\nonumber \\
&&\quad+\,p \,  \rho(i-1,t)\,(1-\rho(i,t))\,
+\,q\, \rho(i+1,t)\,(1-\rho(i,t))
  \nonumber
\end{eqnarray}
with initial condition $\rho(i,0) = \rho_2I(i\leq 0) + \rho_1 I (i >0)$.
The
limiting motion of the rod will be a non-homogeneous (in time) Markov
process described by
\begin{eqnarray}
{d \E (f(Y_t)\,|\,{\cal F}_t)\over dt}&=&
a\,[1-\rho(Y_t+1,t)]^N\,[f(Y_t+1)-f(Y_t)]\nonumber\\
&&\qquad\qquad +\, b\,[1-\rho(Y_t-1,t)]^N\,[f(Y_t-1)-f(Y_t)].\nonumber
\end{eqnarray}
where ${\cal F}_t$ is the sigma field generated by $\{Y_s:s\le
t\}$. These results can be obtained with the techniques we
used to prove Theorem \ref{theo} and will be the content of a future
publication, cf. \cite{pablo}.

\noindent {\bf Remark 2:}  One may wonder
how general the results are.  As an illustration of this we consider
the following somewhat abstract modification of Lemma \ref{conv}.
Suppose that $\mu$ is a reversible measure on $\{0,1\}^{\Z^2}$ both
for a monomer dynamics with generator $\L_1$ and one with generator
$\L_2$. As an example, we could keep in mind the case where $\L_{12} =
\L_1 + \L_2$ is a Kawasaki dynamics (exclusion process with speed
change) at finite temperature with $\L_1$ generating the horizontal
and $\L_2$ generating the vertical jumps; $\mu$ is the corresponding
Gibbs measure. The measure $\mu$ is then also reversible for
$\L_{12}^\gamma = \gamma \L_1 + \L_2$.  Now we insert the polymer and
we fix it at some position $y\in \Z$. The dynamics of the monomers is
now conditioned on having no monomers in the excluded volume $A_N(y)$:
$\eta_t(i)=0, \forall i \in A_N(y), \forall t\geq 0$. The new
generator is $\L_{12}^{y,\gamma} = \gamma \L_1^y + \L_2^y$ obtained by
setting all of the original rates equal to zero for all updating that
would create a monomer in the region $A_N(y)$ (the direct analogue of
what was done in (\ref{hori}) and (\ref{vert})). It follows
then from Lemma \ref{reversi} that
$\mu^y=\mu(\cdot|\eta(i)=0, \forall i \in A_N(y))$ is reversible for
$\L_{12}^y$.  We finally denote by $\mu_x^y$ the restriction of
$\mu^y$ to the layer at height $x$ (i.e., the set $\{i\in \Z^2,
i_2=x\}$). This measure is reversible for $\L_1^{y}$. We have the
following result:
\begin{prop}
\label{4.1}
  Denote by $S^\gamma_y (t)$ the semigroup with generator
  $\L_{12}^{y,\gamma}$. Assume that for all $x$ $\mu_x^y$ is ergodic
  for $\L_1^{y}$.  Let $f_x$ be a function in $L^2{(\mu^y )}$ with
  dependence set on layer $x\not= y$. We have:
\begin{equation}
\lim_{\gamma\uparrow\infty}
\| S^\gamma_y (t) f_x - \int d\mu^y_x\  f_x \|_{L^2 (\mu^y )} = 0.
\end{equation}
\end{prop}
{\bf Proof:}
By ergodicity $\L^y_1$ has simple eigenvalue $0$ with corresponding
eigenspace the constant functions. Hence by the spectral theorem,
\begin{equation}\label{ergospec}
(\int d\mu^y_x f_x )^2 = \E^{\L_{1}^{y}}_{f_x,f_x} ( \{ 0 \} ),
\end{equation}
where $\E^{\L_{1}^{y}}_{f_x,f_x}$ denotes the spectral measure
of the selfadjoint operator $\L_{1}^{y}$. Therefore, we have to show
that if $f_x$ is a function on layer $x$ such that
\begin{equation}\label{nul}
\E^{\L_{1}^{y}}_{f_x,f_x} ( \{ 0 \} )=0,
\end{equation}
then
\begin{equation}\label{semig}
\lim_{\gamma\uparrow\infty} \| S^\gamma_y (t) \|_{L^2 (\mu^y )} = 0.
\end{equation}
For every $\varphi$ in the domain of $\L_{12}^{y,\gamma}$,
\begin{equation}
\lim_{\gamma\uparrow\infty} \frac{1}{\gamma} (\L_{12}^{y,\gamma} \varphi)
= \L_1^y \varphi.
\end{equation}
Hence the spectral measures $\E^{-\frac{1}{\gamma}\L_{12}^{y,\gamma}}_{f_x,f_x}$
converges weakly to the spectral measure
$\E^{-\L^y_1}_{f_x,f_x}$. Therefore, we can estimate
\begin{eqnarray}
\| S^{\gamma}_y (t) f_x \|_{L^2 (\mu^y )}^2 &=&
\int_0^\infty e^{-\gamma t \lambda}
\E^{-\frac{1}{\gamma}\L_{12}^{y,\gamma}}_{f_x,f_x} (d\lambda)\nonumber\\
&\leq &
\int_0^\delta e^{-\gamma t\lambda}
\E^{-\frac{1}{\gamma}\L_{12}^{y,\gamma}}_{f_x,f_x} (d\lambda)
+ e^{-\gamma t\delta} \| f_x\|_{L^2 (\mu^y )}\nonumber\\
&\leq &
\E^{-\frac{1}{\gamma}\L_{12}^{y,\gamma}}_{f_x,f_x}([0,\delta]) +
e^{-\gamma t\delta} \| f_x\|_{L^2 (\mu^y )}.
\end{eqnarray}
Letting $\gamma$ tend to infinity, and then $\delta$ to zero,
using (\ref{nul}), we obtain (\ref{semig}).\QED

\medskip
\noindent {\bf Remark 3:} Proposition \ref{4.1} is general but has a
strong hypothesis: the ergodicity of the one-layer horizontal
dynamics. This is known only in a few cases, in particular in the
symmetric simple exclusion process we treated in Lemma \ref{conv}.
It is however expected to be true at least for high temperature
Kawasaki dynamics.

\end{document}